\documentclass[opre,nonblindrev]{informs3arxiv}
\usepackage{endnotes}
\let\footnote=\endnote
\let\enotesize=\normalsize
\def\notesname{Endnotes}%
\def\makeenmark{$^{\theenmark}$}
\def\enoteformat{\rightskip0pt\leftskip0pt\parindent=1.75em
  \leavevmode\llap{\theenmark.\enskip}}
\usepackage{natbib}
 \bibpunct[, ]{(}{)}{,}{a}{}{,}%
 \def\bibfont{\small}%
\TheoremsNumberedThrough     
\ECRepeatTheorems
\EquationsNumberedThrough    

\usepackage{algorithm}
\usepackage[noend]{algpseudocode}
\usepackage{amsmath}
\usepackage{amssymb}
\usepackage{bm}
\usepackage{algorithm}
\usepackage[noend]{algpseudocode}
\usepackage{booktabs}
\usepackage{tikz}
\usetikzlibrary{shapes,arrows,calc}
\usepackage{csvsimple}
\usepackage{graphicx}
\usepackage{subcaption}
\usepackage{todonotes}
\usepackage{mathtools}
\usepackage{multirow, multicol}

\begin{document}

\RUNAUTHOR{Xavier, Feng and Dey}
\RUNTITLE{Decomposable Formulation of Transmission Constraints}
\TITLE{Decomposable Formulation of Transmission Constraints for Decentralized Power Systems Optimization}

\ARTICLEAUTHORS{%
\AUTHOR{\'Alinson S. Xavier, Feng Qiu}
\AFF{Energy Systems Division, Argonne National Laboratory, Argonne, IL, USA. \EMAIL{axavier@anl.gov}, \EMAIL{fqiu@anl.gov}} 
\AUTHOR{Santanu S. Dey}
\AFF{School of Industrial and Systems Engineering, Georgia Institute of Technology, Atlanta, GA, USA. \EMAIL{santanu.dey@isye.gatech.edu}}
}

\ABSTRACT{%
    One of the most complicating factors in decentralized optimization for power systems is the modeling of power flow equations.
    Existing formulations for direct current (or DC) power flows either have limited scalability or are very dense and unstructured, making them unsuitable for large-scale decentralized studies.
    In this work, we present a novel DC power flow formulation, based on sparsified injection shift factors, which has a decomposable block-diagonal structure, scales well for large systems, and can efficiently handle N-1 security requirements.
    Benchmarks on Multi-Zonal Security-Constrained Unit Commitment problems show that the proposed formulation can reliably and efficiently solve instances with up to 6,515 buses, with no convergence or numerical issues.
}%


\KEYWORDS{Power systems optimization, alternating-method of multipliers (ADMM), injection shift factors, decentralized optimization.}

\maketitle

%


\section{Introduction}

The power industry still predominantly relies on centralized methods for the optimal operation of the electrical grid.
When clearing the day-ahead electricity markets, for example, Independent System Operators (ISOs) collect, in a central location, information coming from a variety of market participants, located across a vast geographical region.
In recent years, such centralized approaches have raised a number of concerns.
First, due to scalability, as the size and complexity of power systems expands, so has the computational power required to optimize their usage. The Midcontinent Independent System Operator (MISO), for example, reported that a recent expansion of market participation and virtual bidding had a serious impact on the computational performance of their software \citep{ChenCastoWangWangWangWan2016}.
The report also describes how simply adding more processing cores to current centralized algorithms is very ineffective at reducing computational times.
The second concern is related to data access. While centralized methods assume that all the necessary data is readily available at one location, the reality is that resources in different regions may belong to different ISOs and Regional Transmission Organizations (RTOs), which may be reluctant to disclose internal details. There is evidence that the lack of computationally efficient methods which can handle data belonging to multiple parties has led to significant resource under-utilization and higher energy prices. A study conducted by ISO New England \citep{ISONE2011}, for instance, showed that, during roughly half of the time in 2009, the power flows in the transmission lines connecting ISO-NE to its neighbors were flowing in the wrong direction.
Decentralized power systems optimization is an alternative approach where, instead of collecting private data from numerous zones into one particular location, each zone independently solves a smaller-scale optimization problem. Some small amount of information is then shared among the zones for coordination, and the subproblems are repeatedly reoptimized, until a global equilibrium is reached. Decentralized optimization can also help other coordination between electricity markets, such as coordinated dispatch, congestion relief, energy exchange, seams study, etc \citep{BALDICK2014319,multiareainterchange}. In contrast to centralized methods, decentralized optimization can naturally scale via parallel computing and does not require sharing of private information between zones.
Over the past decades, various decentralized approaches have been proposed for DC-OPF and Economic Dispatch \citep{wang2016fully,doostizadeh2016multi,loukarakis2015decentralized,phan2014minimal}, AC Optimal Power Flow \citep{erseghe2014distributed,sun2013fully,dall2013distributed,magnusson2015distributed,loukarakis2014investigation} and Unit Commitment \citep{feizollahi2015large,chung2011multi,li2015decentralized}.
We refer to \cite{wang2017distributed} for a more complete survey.

One of the most complicating factors of decentralized optimization for power systems is the modeling of power flow equations. Unlike flows in other domains, electrical power flows have a very complex behavior, dictated by non-linear physical laws. In this manuscript, we focus on optimization problems where the DC linearization of the power flow equations is typically employed, such as Economic Dispatch (ED), Unit Commitment (UC) and Transmission Expansion Planning (TEP). Even when linearized, power flows are still very challenging for decentralized methods, since they have a global effect --- power injections at any particular location may affect flows across the entire network.

Two main formulations for DC power flows have been described in the literature. The \emph{phase-angle formulation} computes the flow in each transmission line based on the phase angle difference between the two buses that the line connects.
Because the flow in each line can be computed from local information only (the phase angles of its endpoints), this is the formulation employed in the vast majority of decentralized methods proposed in the literature. Nevertheless, the phase-angle formulation suffers from two serious limitations which makes it unsuitable for large-scale systems. First, it is unable to exploit the fact that, in realistic systems, only a small fraction of transmission constraints need to be enforced. Second, and more importantly, if the problem requires that flows remain within the safety limits even after the unexpected failure of any individual transmission line --- a requirement known as \emph{N-1 security}, and imposed by the North American Electricity Reliability Corporation (NERC) --- then multiple copies of the phase-angle variables are necessary, making the problem size prohibitively large. For this reason, almost the entire literature on decentralized optimization for power systems has neglected contingencies. One notable exception is \cite{BiskasBakirtzis2004}, which incorporates internal contingencies only, by combining injection shift factors with the formulation we present next.

The \emph{injection shift factors (ISF) formulation}, which is currently the most widely used formulation in the industry, directly calculates the flow in each transmission line as a weighted sum of the net injections at each node of the network. 
Since the flow in each line can be computed independently from the remaining ones, if it is known from previous experience that a certain transmission line is not under risk of exceeding its thermal limits, then the associated flow variables and constraints can be dropped from the formulation. This observation allows the formulation to scale to very large systems, and can also be exploited to efficiently model N-1 security constraints. The main drawback of the ISF formulation, which has prevented its usage in decentralized optimization, is its high density and lack of structure: in order to calculate any flows, precise information of the entire network is needed.

In this work, we describe a novel DC power flow formulation for decentralized optimization, based on geographical decomposition, which combines the best properties of the two formulations above. Like the phase-angle formulation, our proposed formulation is well structured, and computes the power flow in each transmission line using only information located either within the zone, or at the zone's boundary. Like the ISF formulation, our proposed formulation scales very well for large systems, since it allows non-critical variables and constraints to be removed.
The formulation is based on the fact, which we prove in Section~\ref{sec:formulation}, that the injection shift factors of external zones are a convex combination of the shift factors at boundary buses. With a small number of linking constraints, we show that it is possible to drop a large number of coefficients from the ISF matrix, obtaining a much more sparse, well-structured and decomposable formulation.  We also show, in Section~\ref{sec:security}, how can this formulation be adapted to efficiently handle N-1 security requirements. Unlike \cite{BiskasBakirtzis2004}, which only considered internal contingencies, our approach can handle unexpected outages that occur either within or outside the zone. To the best of our knowledge, this is the first work where this has been done.

To evaluate the computational efficiency of the proposed formulation, we benchmark it on the Security-Constrained Unit Commitment Problem (SCUC), a challenging NP-hard problem solved daily to clear the day-ahead electricity markets. In Section~\ref{sec:experiments}, we present computational results on a diverse set of realistic, industrial-sized instances ranging from 1,888 to 6,515 buses, split into two geographical zones.
In our experiments, the proposed formulation is able to reliably solve all instances in under 60 minutes of wallclock time, with no convergence or numerical issues, and very small optimality gaps when compared to centralized solutions. In comparison, the decentralized phase-angle formulation failed to solve even the smallest test cases, even when allowed very relaxed convergence and feasibility tolerances. Although we only present computational results for SCUC, we stress that the formulation is applicable to any problem that uses the DC power flow equations.

\section{Preliminaries}

\subsection{Optimal exchange via ADMM}
\label{subsec:admm}

\emph{Alternative-direction method of multipliers} (ADMM) is an algorithm for solving optimization problems in a distributed computing environment. Since its introduction in the 1970s, ADMM has been successfully applied in a number of fields \citep{boyd2011distributed}, including power systems \citep{erseghe2014distributed,magnusson2015distributed,feizollahi2015large}. In this manuscript, we use it to solve the canonical \emph{optimal exchange problem}, given by
\begin{align}
    \text{minimize} \hspace{2em} & \sum_{i=1}^n f_i(x^i), \label{eq:admm} \\
    \text{subject to} \hspace{2em} & \sum_{i=1}^n x^i = 0, \notag
\end{align}
where $x^i \in \mathbb{R}^k$ and $f_i: \mathbb{R}^k \to \mathbb{R} \cup \{ +\infty \}$, for $i=1,\ldots,n$. The problem is composed by $n$ subsystems, each trying to minimize its own objective function $f_i$, while being restricted by a global equilibrium constraint. When the evaluation of each $f_i$ is computationally expensive, or when the description of these functions relies on private data, it may not be possible to run the entire minimization in a central process. ADMM can then be used, which, in this particular case, reduces to the iterative method described in Algorithm~\ref{alg:admm}. In Step 3 of the algorithm, each subsystem minimizes its own objective function independently and in parallel. In Step 4, a collective all-reduce operation is performed, and all subsystems receive a global average of the private $x^i$ variables. If each component of $x$ is interpreted as the amount of some resource being produced or consumed, then the $\lambda$ variables indicate the prices of these resources. In Step 6, the prices are updated, according to the amount being over- or under-produced. The algorithm stops when a certain convergence tolerance $\varepsilon$ is reached. We note that Algorithm~\ref{alg:admm} can be implemented without centralized coordination. Asynchronous versions of the method have also been proposed. For a more details we refer to \cite[Subsection 7.3.2]{boyd2011distributed}.
In contrast to the \emph{consensus problem}, which has been commonly used for decentralized optimization in power systems, the optimal exchange problem better fits the import/export nature of power exchanges, and, as a consequence, gives more natural and meaningful values to the $\lambda$ variables. We conclude this subsection by recalling that ADMM is only guaranteed to converge when the functions $f_i$ are closed, proper and convex. In other situations, the method can still be used as a heuristic.

\begin{algorithm}
    \begin{algorithmic}[1]
    \normalsize
    \State Let $\rho \in \mathbb{R}$ and $\varepsilon \in \mathbb{R}$ be given.
    \State Let $\tau \leftarrow \textbf{0}_k,
                \lambda \leftarrow \textbf{0}_k,
                \bar{x}^i \leftarrow \textbf{0}_k$,
                for $i=1,\ldots,n$
    \State Let $\bar{x}^i \leftarrow \arg\min_{x^i} \left[
            f_i(x^i) +
            \lambda^T x^i +
            \frac{\rho}{2} \| x^i - \bar{x}^i + \tau \|^2_2 \right] $,
            for $i=1,\ldots,n$
    \State Let $\tau = \frac{1}{n} \sum_{i=1}^n \bar{x}^i$
    \State If $\| \tau \| < \varepsilon$ then stop.
    \State Let $\lambda \leftarrow \lambda + \rho \tau$
    \State Go to step 3.
    \end{algorithmic}
    \caption{Optimal Exchange via ADMM \label{alg:admm}}
\end{algorithm}

\subsection{Centralized phase-angle formulation}
\label{subsec:theta}

Consider a transmission network composed by a set $B$ of buses and a set $L$ of transmission lines. In this manuscript, we represent this network as a directed graph, where the direction assigned to each line is arbitrary. For each bus $b \in B$, we have two decision variables: $n_b$, the net amount of real power (in MW) injected at the bus; and $\theta_b$, the phase angle (in per-unit) at the bus. The \emph{phase-angle formulation} of the DC power flow equations is given by
\begin{subequations}
\begin{align}
        & f_{uv} = \beta_{uv} \left( \theta_u - \theta_v \right) & \forall (u,v) \in L, \label{eq:theta-diff} \\
        & \sum_{u : (u,b) \in L} f_{ub} - \sum_{u : (b,u) \in L} f_{bu} + n_b = 0 & \forall b \in B, \label{eq:theta-flow} \\
        & -F_l \leq f_{l} \leq F_l & \forall l \in L,
\end{align}
\end{subequations}
where $\beta_l$ is the susceptance of transmission line $l$ and $F_l$ is the transmission thermal limit. Constraints \eqref{eq:theta-diff} computes the flow $f_{uv}$ in transmission line $(u,v)$, based on the phase-angle difference of its endpoints. Constraints \eqref{eq:theta-flow} enforce the preservation of flow across the entire network.

\subsection{Decentralized phase-angle formulation}
\label{subsec:dectheta}

We now consider transmission networks decomposed into multiple zones. Consider a partition $(L_1,L_2)$ of the transmission lines $L$ such that the subnetworks induced by $L_1$ and $L_2$ are connected components (that is, they have no islands). We will refer to these subnetworks as zones 1 and 2. Let $B_1$ and $B_2$ be the sets of buses incident only to transmission lines in $L_1$ and $L_2$, respectively. Let $B_\cap$ be set of buses incident to both. We assume that no generators or loads are located at buses $B_\cap$, and therefore the net injection at these buses is zero. This assumption is not restrictive, since any generators or loads located at these buses can be moved to new artificial buses, located either in $B_1$ or $B_2$ and connected to the original bus by artificial transmission lines with very high capacity. The centralized phase-angle formulation, presented in the previous subsection, can be rewritten as
\begin{subequations}
\begin{align}
        & f^k_{uv} = \beta_{uv} \left( \theta^k_u - \theta^k_v \right) & \forall k \in \{1, 2\}, (u,v) \in L_k, \label{eq:disttheta-1} \\
        & \sum_{u : (b,u) \in L_k} f^k_{ub} + \sum_{u : (b,u) \in L_k} f^k_{bu} + n^k_b = 0 & \forall k \in \{1,2\}, b \in B_k, \label{eq:disttheta-2}  \\
        & \sum_{u : (b,u) \in L_k} f^k_{ub} + \sum_{u : (b,u) \in L_k} f^k_{bu} + w^k_b = 0 & \forall k \in \{1,2\}, b \in B_\cap \label{eq:disttheta-3} \\
        & -F_l \leq f^k_{l} \leq F_l
            & \forall k \in \{1,2\}, l \in L_k, \label{eq:disttheta-4} \\
        & w^1_b + w^2_b = 0 & \forall b \in B_\cap \label{eq:disttheta-5} \\
        & \theta^1_b - \theta^2_b = 0& \forall b \in B_\cap \label{eq:disttheta-6} 
\end{align}
\end{subequations}
where $w^k_b$ indicates the amount of power (in MW) exported from (or imported into) zone $k$ at bus $b$.
The superscript $k$ in the formulation above indicates which subsystem, in a distributed computation environment, would own the decision variable. Equations \eqref{eq:disttheta-1}--\eqref{eq:disttheta-4} are clearly local, since all decision variables appearing in these equations belong to the same subsystem. Only equations \eqref{eq:disttheta-5}--\eqref{eq:disttheta-6} affect multiple subsystems.
These equations can be easily handled by the optimal exchange algorithm presented in Subsection~\ref{subsec:admm}.

\subsection{Injection shift factors formulation}
\label{subsec:isf}

The \emph{injection shift factors formulation} of the DC power flow equations is given by
\begin{align*}
    & f_{l} = \sum_{b \in B} \delta_{lb} n_b & \forall l \in L, \\
    & \sum_{b \in B} n_b = 0, \\
    & -F_l \leq f_{l} \leq F_l & \forall l \in L,
\end{align*}
where $\delta_{lb}$ is constant known either as \emph{injection shift factor (ISF)} or \emph{power transfer distribution factor (PTDF)}, and represents the amount of power that flows through line $l$ when 1 MW is injected at $b$ and withdrawn from the slack bus. Let $\Delta \in \mathbb{R}^{|L| \times |B|}$ be the matrix formed by the $\delta_{lb}$ constants. 
The main drawback of the ISF formulation is that $\Delta$ is typically very dense and unstructured. To increase the sparsity of $\Delta$, a common practice in the industry is to discard all $\delta_{lb}$ entries that have magnitude below a fixed threshold. Although this technique can effectively reduce computational times, it results in decreased accuracy and does not improve the structure of $\Delta$.

\section{Decentralized injection shift factors formulation}
\label{sec:formulation}

In this section, we modify the ISF formulation presented in Subsection~\ref{subsec:isf} to make it more suitable for distributed optimization. Such sparse formulations can also have advantages from the perspective of non-decentralized methods (see for example~\cite{bixby2002solving,walter2014sparsity,amaldi2014coordinated,reid1982sparsity,dey2015approximating,dey2018analysis,dey2018theoretical}).

We start in Subsection~\ref{subsec:deccomp} by showing how to compute injection shift factors in a decentralized way. Then, in Subsection~\ref{subsec:decproof}, we show that a large number of coefficients can be dropped from the $\Delta$ matrix by introducing a small number of linking constraints. In Subsection~\ref{subsec:decisf} we present our proposed formulation.

\subsection{Decentralized computation of ISF}
\label{subsec:deccomp}

The first challenge of using the ISF formulation for decentralized optimization is the computation of the injection shift factors. In a centralized setting, the matrix $\Delta$ is typically computed through the expression
\[
    \Delta :=  D M \left(M^T D M\right)^{-1},
\]
where $M \in \{-1,0,1\}^{|L| \times (|B|-1)}$ is incidence matrix of the network, with the column corresponding to the slack bus removed, and $D \in \mathbb{R}^{|L| \times |L|}$ is a diagonal matrix containing the line susceptances. In a decentralized setting, the matrices $D$ and $M$ are not available, and therefore this expression cannot be used.

For decentralized optimization, we propose to compute $\Delta$ by solving $|B|$ Single-Period DC Power Flow problems, using the decentralized phase-angle formulation presented in Subsection~\ref{subsec:dectheta}. In each subproblem, the net injection at the slack bus is set to $-1$ MW, and the net injection of exactly one bus $b \in B$ is set to 1 MW. The vector of flows obtained correspond to one column of the $\Delta$ matrix, by definition. Because these subproblems are single-period and have an empty objective function, they can be solved very efficiently in practice, even using the phase-angle formulation. Since the subproblems have no interdependencies, they can also be solved in parallel, making this task even faster.

We also clarify that $\Delta$ only needs to be computed once for each transmission network. Once this matrix is pre-computed through the simplified problem above, it can be repeatedly used to formulate much more challenging power systems optimization problems, such as Multi-Period DC Optimal Power Flow (DC OPF) and Security-Constrained Unit Commitment (SCUC).

\subsection{Sparsifying the ISF matrix}
\label{subsec:decproof}

Even if the matrix of injection shift factors $\Delta$ is available, it still lacks the block-diagonal structure that is typically required for decentralized methods to work well. In this subsection we prove that, when the transmission network can be partitioned into zones that have only a small number of tie lines between them, then, by adding a small number of linking constraints to the original ISF formulation, a large number of coefficients in the $\Delta$ matrix can be dropped, resulting in a much more sparse and well-structured matrix. For simplicity, we focus on the 2-zone case. Each zone can be further subdivided, if desired, using the same method proposed in this subsection.

Let $(L_1,L_2)$ and $(B_1,B_2,B_\cap)$ be a partition of the transmission network as defined as in Subsection~\ref{subsec:dectheta}. In the following, we assume that the rows of $\Delta$ corresponding to lines $L_1$ were computed with the slack bus located in $B_1$, while the rows corresponding to lines $L_2$ were computed with a slack bus in $B_2$. Since, in the ISF formulation presented in Subsection~\ref{subsec:isf}, each transmission constraint is completely independent from the remaining ones, this usage of multiple slack buses is allowed, and does not change the values of the $f_l$ variables.
For clarity, we partition the $\Delta$ matrix as:
\[
    \Delta = \left[
        \begin{array}{c|c|c}
        \hspace{0.2in} \Delta^{11}    \hspace{0.2in} &
        \hspace{0.1in} \Delta^{1\cap} \hspace{0.1in} &
        \hspace{0.2in} \Delta^{12}    \hspace{0.2in} \\
        \hline
        \Delta^{21} & \Delta^{2\cap} & \Delta^{22}
        \end{array}
    \right]
\]

\newcommand{\red}{\bf}
 \begin{figure}
    \caption{Numerical example of partitioned Injection Shift Factor (ISF) matrix for
    a network with 8 buses. \label{fig:ex}}
    \begin{subfigure}[c]{0.7\textwidth}
        \centering
        \setlength{\tabcolsep}{0.5em}
        \scriptsize
        \begin{tabular}{l|rrr|rr|rrr}
        \toprule
        {} &  1 &     2 &     3 &     4 &     5 &     6 &     7 &     8 \\
        \midrule
        (1, 2) &  0 & -0.65 & -0.06 & \red -0.29 & \red -0.12 & -0.23 & \red -0.19 & -0.21 \\
        (1, 3) &  0 & -0.06 & -0.81 & \red -0.12 & \red -0.62 & -0.31 & \red -0.44 & -0.38 \\
        (1, 4) &  0 & -0.29 & -0.12 & \red -0.58 & \red -0.25 & -0.46 & \red -0.38 & -0.42 \\
        (2, 4) &  0 &  0.35 & -0.06 & \red -0.29 & \red -0.12 & -0.23 & \red -0.19 & -0.21 \\
        (3, 5) &  0 & -0.06 &  0.19 & \red -0.12 & \red -0.62 & -0.31 & \red -0.44 & -0.38 \\
        \hline
        (4, 6) &  0.62 &  0.69 &  0.44 &  0.75 &  0.25 & -0.06 &  0.06 &  0 \\
        (5, 7) &  0.38 &  0.31 &  0.56 &  0.25 &  0.75 &  0.06 & -0.06 &  0 \\
        (6, 7) &  0.08 &  0.12 & -0.04 &  0.17 & -0.17 &  0.29 & -0.29 &  0 \\
        (6, 8) &  0.54 &  0.56 &  0.48 &  0.58 &  0.42 &  0.65 &  0.35 &  0 \\
        (7, 8) &  0.46 &  0.44 &  0.52 &  0.42 &  0.58 &  0.35 &  0.65 &  0 \\
        \bottomrule
        \end{tabular}
        \caption{Partitioned ISF. \label{fig:ex-isf}}
    \end{subfigure}
    \begin{subfigure}[c]{0.25\textwidth}
        \centering
        \begin{tikzpicture}
            \pgfsetlinewidth{0.8pt}
            \scriptsize
        
            \coordinate[label={[label distance=3pt]135:1}] (b1) at (0, 1);
            \coordinate[label={[label distance=3pt]135:2}] (b2) at (1, 2);
            \coordinate[label={[label distance=3pt]200:3}] (b3) at (1, 0);
            \coordinate[label={[label distance=3pt]95:4}] (b4) at (2, 1.5);
            \coordinate[label={[label distance=3pt]135:5}] (b5) at (2, 0.5);
            \coordinate[label={[label distance=3pt]135:6}] (b6) at (3, 2);
            \coordinate[label={[label distance=3pt]200:7}] (b7) at (3, 0);
            \coordinate[label={[label distance=3pt]45:8}] (b8) at (4, 1);
        
            \draw (b1) -- (b2);
            \draw (b1) -- (b3);
            \draw (b1) -- (b4);
            \draw (b2) -- (b4);
            \draw (b3) -- (b5);
            \draw (b4) -- (b6);
            \draw (b5) -- (b7);
            \draw (b6) -- (b7);
            \draw (b6) -- (b8);
            \draw (b7) -- (b8);
        
            \draw[dashed, black] (2, -1) -- (2, 3);
        
            \draw[fill=black] (b1) circle[radius=0.1];
            \draw[fill=white] (b2) circle[radius=0.1];
            \draw[fill=white] (b3) circle[radius=0.1];
            \draw[fill=white] (b4) circle[radius=0.1];
            \draw[fill=white] (b5) circle[radius=0.1];
            \draw[fill=white] (b6) circle[radius=0.1];
            \draw[fill=white] (b7) circle[radius=0.1];
            \draw[fill=black] (b8) circle[radius=0.1];
        \end{tikzpicture}
        \caption{Simplified network representation. \label{fig:ex-graph}}
    \end{subfigure}
\end{figure}

In the following, we show that each column of $\Delta^{12}$ and $\Delta^{21}$ is a convex combination of the columns of $\Delta^{1\cap}$ and $\Delta^{2\cap}$, respectively. First, we present a numerical example to clarify this result.

\begin{example}
Figure~\ref{fig:ex-graph} shows a simplified representation of a transmission network with eight buses and two zones. Each bus is represented as a vertex in a graph, and each transmission line is represented as an edge. Slack buses are represented in black, and, for demonstration purposes, all transmission lines have the same susceptance, although this is not a required assumption. In this example, $B_1=\{1,2,3\}, B_\cap=\{4,5\}$ and $B_2=\{6,7,8\}$. Figure~\ref{fig:ex-isf} shows the partitioned ISF matrix $\Delta$ corresponding to this network. In the following, we prove that the upper half of column 7, highlighted in bold, is a convex combination of the upper halves of columns 4 and 5. More precisely,
\[
    0.38 \left[ \begin{matrix} -0.29 \\ -0.12 \\ -0.58 \\ -0.29 \\ -0.12 \end{matrix}\right]
    + 0.62 \left[ \begin{matrix} -0.12 \\ -0.62 \\ -0.25 \\ -0.12 \\ -0.62 \end{matrix} \right]
    \approx
    \left[ \begin{matrix} -0.19 \\ -0.44 \\ -0.38 \\ -0.19 \\ -0.44 \end{matrix} \right].
\]
Here, we use approximate equality because all entries of the matrix in Figure~\ref{fig:ex-isf} have been rounded to two decimals digits; if the matrix is represented exactly, equality holds.
\hfill $\diamond$
\end{example}

\begin{theorem}
    \label{thm:decomposition}
    Let $k \in \{1,2\}$ and $c \in B_k$. There exist constants $\gamma_{bc}^k \in [0, 1]$, for $b \in B_\cap$ such that:
    \begin{subequations}
    \begin{align}
        & \delta_{lc} = \sum_{b \in B_\cap} \delta_{lb} \gamma^k_{bc}
            & \forall l \in L_{3-k}, \label{eq:thm-a} \\
        & \sum_{b \in B_\cap} \gamma_{bc}^k = 1 \label{eq:thm-b}
    \end{align}
    \end{subequations}
\end{theorem}
\proof{Proof.}
    For simplicity, we assume $k=2$. Let $s \in B_1$ be the slack bus from zone 1. We also define $\hat{B} = B \setminus \{s\}$ and $\hat{B}_1 = B_1 \setminus \{s\}$. Let $M \in \{-1,0,1\}^{L \times \hat{B}}$ be the (reduced) incidence matrix of the network, let $D \in \mathbb{R}^{L \times L}$ be the diagonal matrix of line susceptances and let $G = M^T D M$. Furthermore, let $n \in \mathbb{R}^{\hat{B}}$ be a column vector such that, for every $b \in \hat{B}$, we have
    \[
        n_b = \begin{cases}
            1 & \text{if } b = c, \\
            0 & \text{otherwise.}
        \end{cases}
    \]
    By definition, $\delta_{\bullet b}$ is the vector of line flows when the net injection at bus $s$ is $-1$ MW, and the net injections of the remaining buses is $n$. Let $\theta \in \mathbb{R}^{\hat{B}}$ be the vector of phase angles in this scenario. That is, let $\theta \in \mathbb{R}^{\hat{B}}$ be such that
    \begin{equation}
        \label{eq:proof-a}
        G \left[ \begin{matrix}
            \theta_{\hat{B}_1} \\
            \theta_{B_\cap} \\
            \theta_{B_2} \\
        \end{matrix} \right] = n.
    \end{equation}
    The next claim shows that, by modifying only the phase angles in $B_2$, it is possible to shift the 1 MW net injection from bus $c \in B_2$ to buses the buses in $B_\cap$.
    \begin{claim}
        There exist $\tilde{\theta} \in \mathbb{R}^{\hat{B}}$ and $\tilde{n} \in \mathbb{R}^{\hat{B}}$ such that $\tilde{\theta}_{\hat{B}_1} = \theta_{\hat{B}_1}, \tilde{\theta}_{B_\cap} = \theta_{B_\cap}, \tilde{n}_{\hat{B}_1} = \bm{0}, \tilde{n}_{B_2} = \bm{0},$ and
        \begin{equation}
            \label{eq:proof-b}
            G \tilde{\theta} =
            G \left[ \begin{matrix}
                \theta_{\hat{B}_1} \\
                \theta_{B_\cap} \\
                \tilde{\theta}_{B_2} \\
            \end{matrix} \right] = \left[ \begin{matrix}
                \bm{0} \\
                \tilde{n}_{B_\cap} \\
                \bm{0}
            \end{matrix} \right] \coloneqq
            \tilde{n}.
        \end{equation}
    \end{claim}
    \proof{Proof of the claim.}
        Partitioning the rows and columns of $G$, we may rewrite \eqref{eq:proof-b} as
        \begin{equation}
            \label{eq:proof-c}
            \left[
                \begin{matrix}
                G^{11} & G^{1\cap} & G^{12} \\
                G^{\cap 1} & G^{\cap \cap} & G^{\cap 2} \\
                G^{21} & G^{2\cap} & G^{22}
                \end{matrix}
            \right]
            \left[ \begin{matrix}
                \theta_{\hat{B}_1} \\
                \theta_{B_\cap} \\
                \tilde{\theta}_{B_2} \\
            \end{matrix} \right] = \left[ \begin{matrix}
                \bm{0} \\
                \tilde{n}_{B_\cap} \\
                \bm{0}
            \end{matrix} \right].
        \end{equation}
        We recall that $G$ is a Laplacian matrix with one column and one row (corresponding to the slack bus) removed. Since there are no edges between $B_1$ and $B_2$, all entries of $G^{12}$ and $G^{21}$ are zero. Therefore,
        \[
            G^{11} \theta_{\hat{B}_1} + G^{1\cap} \theta_{B_\cap} + G^{12} \tilde{\theta}_{B_2}
            = G^{11} \theta_{\hat{B}_1} + G^{1\cap} \theta_{B_\cap}
            = \bm{0},
        \]
        where the last equality follows from \eqref{eq:proof-a}. This shows that the first set of constraints in \eqref{eq:proof-c} is always satisfied. Rewriting the two remaining sets of constraints, we obtain
        \begin{equation}
            \label{eq:proof-d}
            \underbrace{\left[ \begin{matrix}
                G^{\cap 2} & -\mathcal{I} \\
                G^{22} & \bm{0}
            \end{matrix} \right]}_{Q}
            \left[ \begin{matrix}
                \tilde{\theta}_{B_2} \\
                \tilde{n}_{B_\cap}
            \end{matrix} \right] =
            -\left[ \begin{matrix}
                G^{\cap 1} & G^{\cap \cap} \\
                \bm{0} & G^{\cap 2}
            \end{matrix} \right]
            \left[ \begin{matrix}
                \theta_{\hat{B}_1} \\
                \theta_{B_\cap} \\
            \end{matrix} \right],
        \end{equation}
        where $\mathcal{I}$ is the identify matrix.
        Since $G^{22}$ is non-singular, then $Q$ is also non-singular. This implies that this system of linear equations always has a solution, proving that the desired $\tilde{\theta}$ and $\tilde{n}$ always exist.
        \hfill$\blacksquare$
    \endproof
    Let $\tilde{\theta}$ and $\tilde{n}$ be vectors satisfying \eqref{eq:proof-b} and let $\gamma^2_{bc} = \tilde{n}_b$, for every $b \in B_\cap$.
    Now we prove that $\gamma^2$ satisfies \eqref{eq:thm-a} and \eqref{eq:thm-b}. Consider a transmission line $l \in L_1$. Recall that $\delta_{lc}$ is the flow in $l$ when the phase angles and net injections are $\theta$ and $n$, respectively. Let $\tilde{\delta}_{lc}$ be the flow in $l$ when the phase angles and net injections are $\tilde{\theta}$ and $\tilde{n}$.
    Computing $\tilde{\delta}_{lc}$ through injection shift factors, we have
    \[
        \tilde{\delta}_{lc}
        = \sum_{b \in \hat{B}} \delta_{lb} \tilde{n}_b
        = \sum_{b \in B_\cap} \delta_{lb} \tilde{n}_b
        = \sum_{b \in \hat{B}} \delta_{lb} \gamma^2_{bc}.
    \]
    Note, however, that the phase angles at the endpoints of $l$ have not changed, and therefore $\delta_{lc} = \tilde{\delta}_{lc}$. This proves that $\gamma^2$ satisfies \eqref{eq:thm-a}. Since the slack bus $s$ is located in $B_1$, the flows in its incident transmission lines have also not changed, and therefore its net injection is still $-1$ MW.
    Flow preservation implies
    \[
        1
        = \sum_{b \in \hat{B}} \tilde{n}_b
        = \sum_{b \in B_\cap} \tilde{n}_b
        = \sum_{b \in B_\cap} \gamma^2_{bc},
    \]
    proving that $\gamma^2$ also satisfies \eqref{eq:thm-b}.
    \hfill\halmos
\endproof

\subsection{Decentralized ISF formulation}
\label{subsec:decisf}

Based on Theorem~\ref{thm:decomposition}, we can rewrite the transmission constraints in the ISF formulation, as well as the power balance equation, in a more decomposable way. For every $b \in B_\cap$, let $w^1_b, w^2_b$ be auxiliary decision variables. As in Subsection~\ref{subsec:dectheta}, we replace the decision variables $n_b$ by either $n^1_b$ or $n^2_b$, depending on whether $b \in B_1$ or $b \in B_2$, to indicate which subsystem, in a distributed computing environment, would own the decision variable. We still make the assumption that $n_b=0$ for every $b \in B_\cap$. Our proposed \emph{decentralized injection shift factor formulation} is given by:
\begin{subequations}
    \label{eq:decisf}
    \begin{align}
        & w^1_b = \sum_{c \in B_2} \gamma^2_{bc} n^2_c
            & \forall b \in B_\cap, \label{eq:decisf-a} \\
        & w^2_b = \sum_{c \in B_1} \gamma^1_{bc} n^1_c
            & \forall b \in B_\cap, \label{eq:decisf-b} \\
        & \sum_{b \in B_k} n^k_b + \sum_{b \in B_\cap} w^k_b = 0
            & \forall k \in \{1,2\}, \label{eq:decisf-c} \\
        & f^k_l =
            \sum_{b \in B_k} \delta_{lb} n^k_b +
            \sum_{b \in B_\cap} \delta_{lb} w^k_b
            & \forall k \in \{1,2\}, l \in L_k \label{eq:decisf-d} \\
        & -F_l \leq f^k_{l} \leq F_l & \forall k \in \{1,2\}, l \in L_k, \label{eq:decisf-e}
    \end{align}
\end{subequations}
By applying Theorem~\ref{thm:decomposition}, it can be easily verified that this formulation is equivalent to the original ISF formulation presented in Subsection~\ref{subsec:isf}. In the proposed formulation, however, assuming that $B_\cap$ is small, only a small number of constraints, namely \eqref{eq:decisf-a} and \eqref{eq:decisf-b}, are non-local. Equations \eqref{eq:decisf-c}--\eqref{eq:decisf-e}, which comprise the vast majority of constraints for large-scale systems, are now completely local. Similarly to the original formulation, note that the $f^k_l$ variables are still independent from each other, and therefore can be added lazily to the formulation, making the proposed formulation scalable to very large systems.

\begin{figure}
    \centering
    \caption{Comparison between sparsity patterns of different formulations. \label{fig:sparsity}}
    \begin{subfigure}[b]{0.4\textwidth}
        \centering
        \fbox{\includegraphics[width=\textwidth]{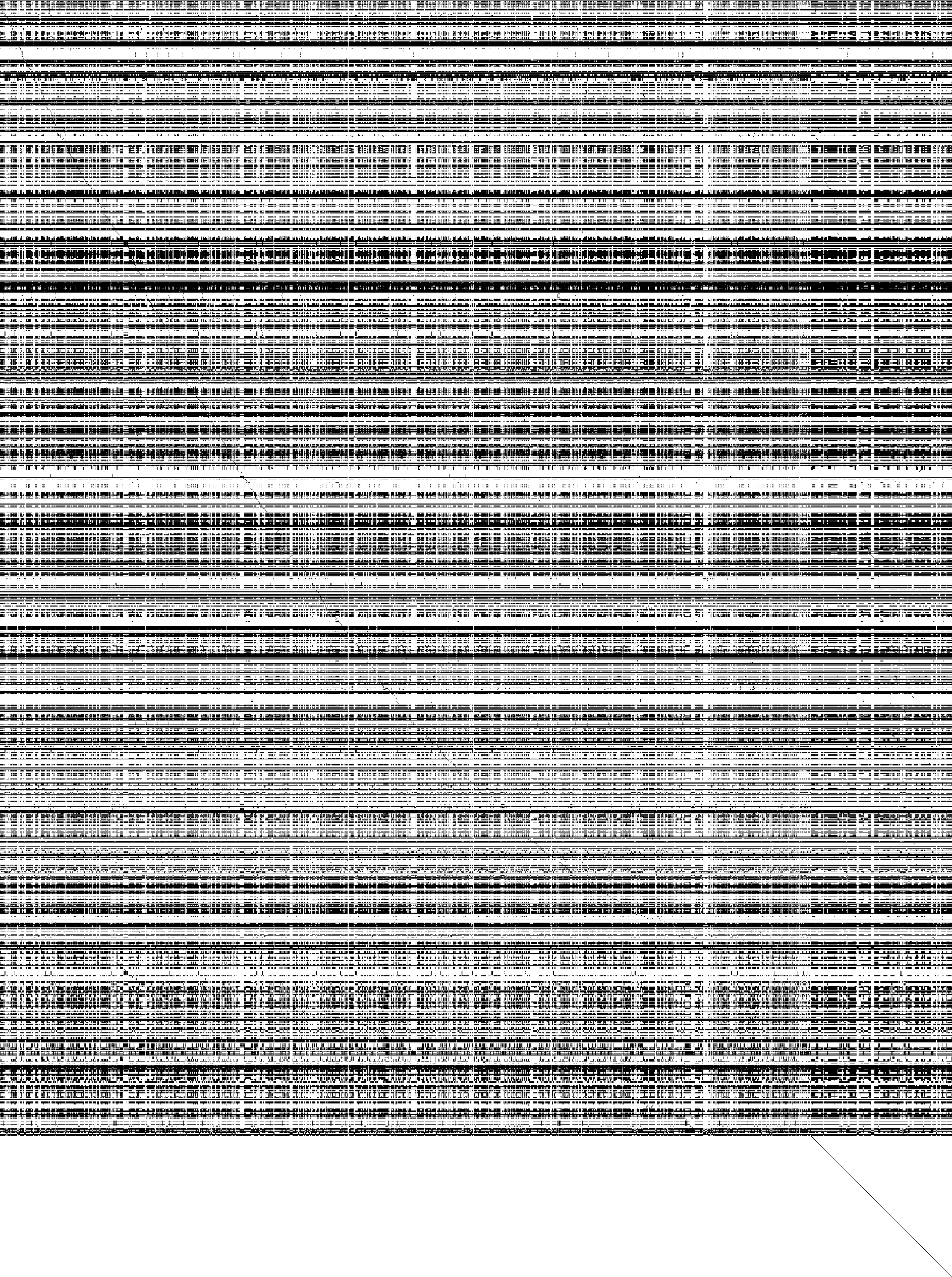}}
        \caption{Original ISF formulation.}
    \end{subfigure}~~~
    \begin{subfigure}[b]{0.4\textwidth}
        \centering
        \fbox{\includegraphics[width=\textwidth]{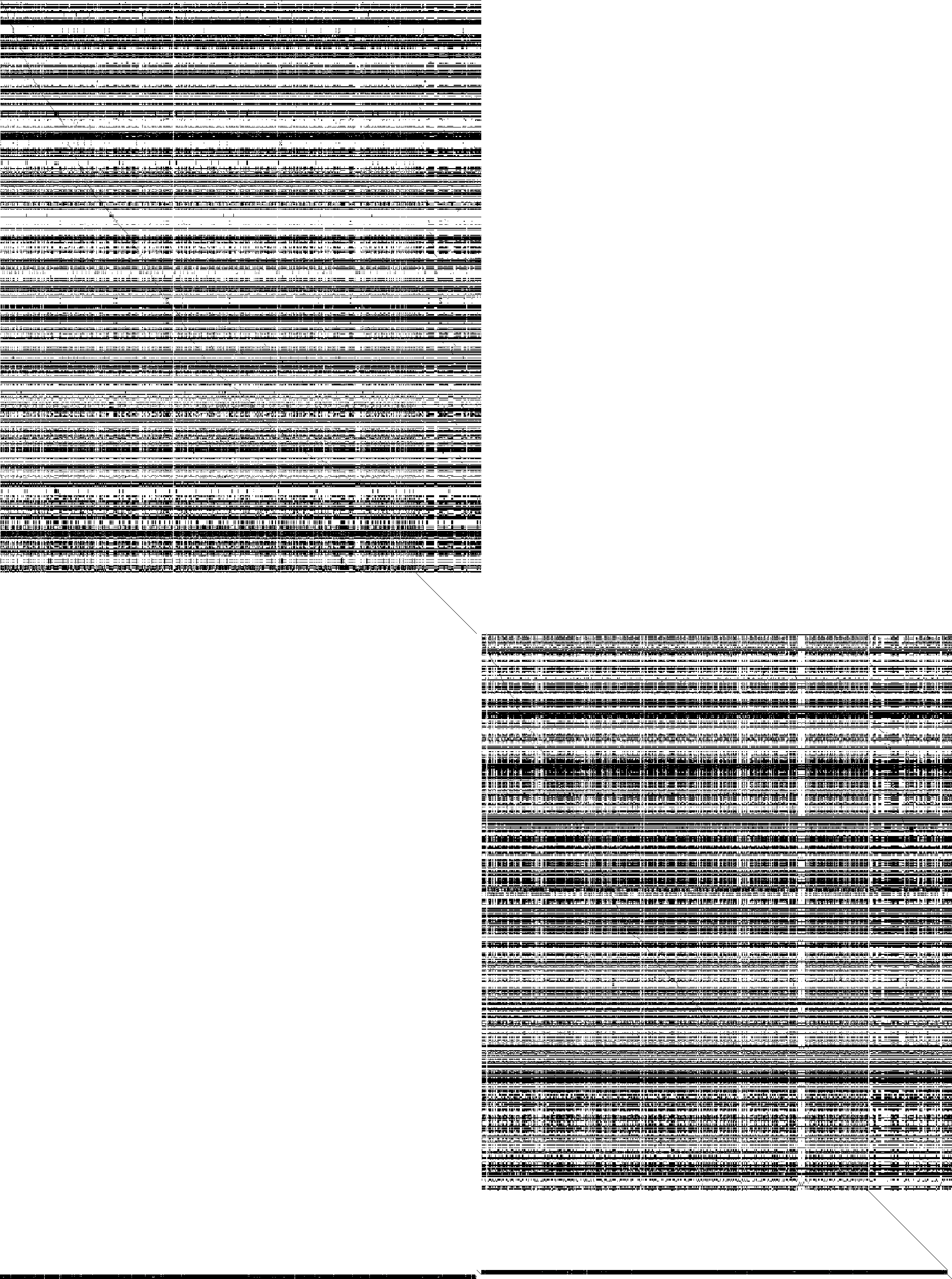}}
        \caption{Proposed formulation. \label{fig:sparsity-b}}
    \end{subfigure}
\end{figure}

Figure~\ref{fig:sparsity} shows a visual comparison between the sparsity structures of the original ISF formulations versus our proposed decomposed version, for the \texttt{case1888rte} instance, which corresponds to the French VHV System in 2013. In the diagram, each non-zero constraint coefficient is represented as a black dot. While the original formulation has no clear structure, the decomposed formulation presentes a clear block-diagonal structure. A very small number of (non-local) linking constraints can be seen at the bottom of Figure~\ref{fig:sparsity-b}.

\section{Enforcement of N-1 security constraints}
\label{sec:security}

In many power system optimization problems, the optimal net injections must induce a feasible network flow not only for the original transmission network, where all transmission lines are operational, but also for a number of N-1 contingency scenarios, where exactly one transmission line has unexpectedly failed.  In this section we show how can the formulation presented in Subsection~\ref{subsec:decisf} be adapted to enforce such constraints.

Consider a pair of transmission lines $m,q \in L$. We will refer to $m$ as the \emph{monitored} line, and $q$ as the \emph{outaged} line. We would like to enforce the transmission limits on $m$, whether $q$ is operational or not. We assume that disconnecting $q$ does not create islands, and that all net injections remain the same. In Subsection~\ref{subsec:int-outage}, we consider the case where $m$ and $q$ belong to the same zone. Then, in Subsection~\ref{subsec:ext-outage}, we consider the case where they belong to different zones.

\subsection{Internal outages}
\label{subsec:int-outage}

Suppose that the monitored line $m$ and the outaged line $q$ belong to the same zone.
Let $f_m$ and $f_q$ be the flows in $m$ and $q$ when $q$ is operational. Also, let $\tilde{f}_{mq}$ be the flow in transmission line $m$ after $q$ is disconnected. It is well known that there exists a constant $\phi_{mq}$, known as \emph{line outage distribution factor}, such that
\[
    \tilde{f}_{mq} = f_m + \phi_{mq} f_q.
\]
The constant $\phi_{mq}$ can be easily computed from the matrix of injection shift factors \citep{guler2007generalized}.
To enforce valid flows under all internal outages, it is sufficient to add the following constraints to Formulation \eqref{eq:decisf}:
\begin{align*}
    & \tilde{f}^k_{mq} = f^k_m + \phi_{mq} f^k_q
        & \forall k \in \{1,2\}, m \in L_k, q \in L_k \\
    & -F_{m} \leq \tilde{f}^k_{mq} \leq F_{m}
        & \forall k \in \{1,2\}, m \in L_k, q \in L_k
\end{align*}
Similar to the regular transmission constraints, if it is known, from previous experience, that a certain flow $\tilde{f}^k_{mq}$ is not under risk of exceeding its limits, this variable and its associated constraints can be dropped from the formulation.

\subsection{External outages}
\label{subsec:ext-outage}

Now consider the case where the monitored transmission line $m$ and the outage transmission line $q$ belong to different zones. In this situation, the approach outlined in the Subsection~\ref{subsec:int-outage} cannot be used to compute $\tilde{f}_{mq}$, since $f_m$ and $f_q$ belong to different subproblems.
This could be solved by requiring both zones to share with each other the precise flows in each of their own transmission lines. This solution, however, does not scale well computationally, since it dramatically increases the number of decision variables that need to reach consensus. It may also not be acceptable for privacy reasons. Because of these difficulties, previous decentralized formulations such as \cite{BiskasBakirtzis2004} simply ignore external outages. In this subsection, we present an alternative solution, with better scalability, which requires significant less sharing.

Let $T$ denote the set of alternative network topologies in which the transmission limits need to be enforced. More precisely, each $\tau \in T$ corresponds to an alternative transmission network containing the same set of buses $B$, but having exactly one transmission line $q \in L$ removed.
For each $\tau \in T$, let $f_{l\tau}$ be the flow in transmission line $l$ under topology $\tau$.
Let $\delta_{lb\tau}$ be the injection shift factors in topology $\tau$, and let $\gamma_{bc\tau}$ be the constants from Theorem~\ref{thm:decomposition} when applied to $\delta_{lbt}$. To enforce transmission limits
in all topologies $T$, one possible solution would be to add the following constraints to Formulation~\eqref{eq:decisf}:
\begin{subequations}
    \label{eq:topology1}
    \begin{align}
        & w^1_{b\tau} = \sum_{c \in B_2} \gamma^2_{bc\tau} n^2_c
            & \forall b \in B_\cap, \tau \in T \label{eq:topology1-a} \\
        & w^2_{b\tau} = \sum_{c \in B_1} \gamma^1_{bc\tau} n^1_c
            & \forall b \in B_\cap, \tau \in T, \label{eq:topology1-b} \\
        & \sum_{b \in B_k} n^k_b + \sum_{b \in B_\cap} w^k_{b\tau} = 0
            & \forall k \in \{1,2\}, \tau \in T \label{eq:topology1-c} \\
        & f^k_{l\tau} =
            \sum_{b \in B_k} \delta_{lb\tau} n^k_b +
            \sum_{b \in B_\cap} \delta_{lb\tau} w^k_{b\tau}
            & \forall k \in \{1,2\}, l \in L_k, \tau \in T \label{eq:topology1-d} \\
        & -F_l \leq f^k_{l\tau} \leq F_l & \forall k \in \{1,2\}, l \in L_k, \tau \in T. \label{eq:topology1-e} 
    \end{align}
\end{subequations}
Although it is possible, in theory, to solve this formulation using the ADMM procedure from Subsection \ref{subsec:admm}, we expect very poor computational performance, due to the large number of consensus variables $w^k_{b\tau}$. In the following, we propose an alternative solution method to enforce these constraints, which does not require any additional consensus variables.
We start by showing that Equations \eqref{eq:topology1-c} are redundant, and therefore can be omitted.
\begin{proposition}
    \label{thm:topology1}
    Let $(w,n,f)$ be a solution to Equations \eqref{eq:decisf}, \eqref{eq:topology1-a} and \eqref{eq:topology1-b}. Then
    \[
        \sum_{b \in B_\cap} w^k_b = \sum_{b \in B_\cap} w^k_{b\tau},
    \]
    for every $k \in \{1,2\}, \tau \in T$.
\end{proposition}
\proof{Proof.}
    Let $\tau \in T$. Without loss of generality, we assume $k=1$. By Theorem~\ref{thm:decomposition},
    \[
        \sum_{b \in B_\cap} \gamma^2_{bc}
        = 1
        = \sum_{b \in B_\cap} \gamma^2_{bc\tau}.
    \]
    Then we have
    \[
        \sum_{b \in B_\cap} w^1_b
        = \sum_{b \in B_\cap} \sum_{c \in B_2} \gamma^2_{bc} n^2_c
        = \sum_{c \in B_2} \left[ n^2_c \sum_{b \in B_\cap} \gamma^2_{bc} \right]
        = \sum_{c \in B_2} \left[ n^2_c \sum_{b \in B_\cap} \gamma^2_{bc\tau} \right]
        = \sum_{b \in B_\cap} w^1_{b\tau}.
    \]
    \hfill\halmos
\endproof
Next, we focus on simplifying Equations \eqref{eq:topology1-d} and \eqref{eq:topology1-e}. For every $b \in B_\cap$ and $\tau \in T$, let $e^1_{b\tau}$ be a new decision variable representing the difference between $w^1_{b\tau}$ and $w^1_b$. That is,
\[
    e^1_{b\tau} := w^1_{b\tau} - w^1_b = \sum_{b \in B_2} \left(
        \gamma^2_{bc\tau} - \gamma^2_{bc}
    \right) n^2_c.
\]
Let $e^2_{l\tau}$ be similarly defined. Additionally, for every $k \in \{1,2\}, l \in L_k$ and $\tau \in T$, let $g^k_{l\tau}$ be a new decision variable representing the difference between $f^k_{l\tau}$ and $f^k_l$. That is,
\begin{align*}
    g^k_{l\tau}
    & := f^k_{l\tau} - f^k_l
    = \sum_{b \in B_k} \left(
        \delta_{lb\tau} - \delta_{lb}
    \right) n^k_b + \sum_{b \in B_\cap} \left(
        \delta_{lb\tau} w^k_{b\tau} - \delta_{lb} w^k_b
    \right) \\
    & = \sum_{b \in B_k} \left(
        \delta_{lb\tau} - \delta_{lb}
    \right) n^k_b + \sum_{b \in B_\cap} \left[
        \delta_{lb\tau} \left( e^k_{b \tau} + w^k_b \right)
        - \delta_{lb} w^k_b
    \right] \\
    & = \sum_{b \in B_k} \left(
        \delta_{lb\tau} - \delta_{lb}
    \right) n^k_b + \sum_{b \in B_\cap} \left[
        \left(
            \delta_{lb\tau} - \delta_{lb}
        \right) w^k_b + \delta_{lb\tau} e^k_{b\tau}
    \right]
\end{align*}
With these auxiliary variables, it can be easily verified that Equations \eqref{eq:topology1-a}, \eqref{eq:topology1-b}, \eqref{eq:topology1-d} and \eqref{eq:topology1-e} can be replaced by
\begin{align}
    & -F_l \leq f^k_{l} + g^k_{l\tau} \leq F_l & \forall k \in \{1,2\}, l \in L_k, \tau \in T.
    \label{eq:topology2}
\end{align}
What we propose next is to replace, in the definitions of $e^k_{l\tau}$ and $g^k_{l\tau}$, the decision variables $n^k_b$ and $w^k_b$ by constants $\tilde{n}^k_b$ and $\tilde{w}^k_b$, corresponding to the the optimal values of $n^k_b$ and $w^k_b$ in the previous ADMM iteration.
During the first ADMM iteration, $\tilde{n}^k_b$ and $\tilde{w}^k_b$ are set to zero. 
With this modification, the variables $e^k_{l\tau}$ and $g^k_{l\tau}$ become constants, which we denote by $\tilde{e}^k_{k\tau}$ and $\tilde{g}^k_{l\tau}$. For clarity, these constants are defined as
\begin{align*}
    & \tilde{e}^1_{b\tau} := \sum_{c \in B_2} \left(
        \gamma^2_{bc\tau} - \gamma^2_{bc}
    \right) \tilde{n}^2_c,
    & \forall b \in B_\cap, \tau \in T
    \\
    & \tilde{e}^2_{b\tau} := \sum_{c \in B_1} \left(
        \gamma^1_{bc\tau} - \gamma^1_{bc}
    \right) \tilde{n}^1_c,
    & \forall b \in B_\cap, \tau \in T
    \\
    & \tilde{g}^k_{l\tau}
    := \sum_{b \in B_k} \left(
        \delta_{lb\tau} - \delta_{lb}
    \right) \tilde{n}^k_b + \sum_{b \in B_\cap} \left[
        \left(
            \delta_{lb\tau} - \delta_{lb}
        \right) \tilde{w}^k_b + \delta_{lb\tau} \tilde{e}^k_{b\tau}
    \right]
    & \forall k \in \{1,2\}, l \in L_k, \tau \in T
\end{align*}
Converting $g^k_{l\tau}$ to a constant causes the vast majority of Equations~\eqref{eq:topology2} to become redundant. Indeed, it is sufficient to enforce only the following set of constraints:
\begin{equation}
    \label{eq:robust}
    -F_l - \underbrace{\left(
            \min_{\tau \in T} \tilde{g}^k_{l\tau}
        \right)}_{\tilde{g}^k_{l,\min}}
        \leq
        f^k_{l}
        \leq F_l - \underbrace{\left(
            \max_{\tau \in T} \tilde{g}^k_{l\tau}
        \right)}_{\tilde{g}^k_{l,\max}}
\end{equation}
for every $k \in \{1,2\}$ and $l \in L_k$.
Equation~\eqref{eq:robust} can be seen as a robust version of Equation \eqref{eq:decisf-e}.

Given these modifications, our proposed solution method is the following. At the end of each ADMM iteration, each zone computes and shares with each other the $\tilde{e}^k_{b\tau}$ values. Unlike network flows, these values reveal very little private information. Upon receiving the updated $\tilde{e}^k_{b\tau}$ values, each zone computes $\tilde{g}^k_{l,\min}$ and $\tilde{g}^k_{l,\max}$ for each transmission line and updates the bounds of Equation~\eqref{eq:robust}. The subproblems are then reoptimized, and the procedure repeats until the solutions converge.

\section{Computational performance}
\label{sec:experiments}

The computational performance of the proposed formulation was evaluated on the Security-Constrained Unit Commitment Problem (SCUC), a challenging NP-hard problem used, for example, to clear the day-ahead electricity markets. Besides transmission and N-1 security, other enforced constraints included (i) maximum and minimum generation limits, (ii) ramping restrictions and (iii) minimum uptime and downtime. 
In our experiments, the decentralized phase-angle formulation from Subsection~\ref{subsec:dectheta} was used as a baseline. Both formulations were implemented in Julia 1.5 and JuMP 0.21, and shared the same ADMM code. IBM ILOG CPLEX 12.9 was used as MIQP and QP solver. Experiments were run a desktop computer (AMD Ryzen 9 3950X, 16 cores, 32 threads, 3.5 GHz, 64 GB DDR4). A single instance was solved at a time, with one process per zone. Inter-process communication was performed via MPI. A wallclock time limit of 3600s was imposed over the entire optimization process.

\subsection{Instances}
Seven instances from MATPOWER \citep{matpower}, corresponding to realistic, large-scale European test systems, were selected to evaluate the formulations. Table~\ref{table:size} presents their main characteristics, including number of buses, generators and transmission lines. Some generator data, such as ramping rates, was missing from the original instances, and was artificially generated based on real data distributions.
The augmented instances have been made available as part of the open-source package UnitCommitment.jl \citep{UCJL}.
To split each instance into two zones, an auxiliary MILP was solved, as described in the appendix. No attempt was made to keep the number of generators in different zones balanced, although the auxiliary problem could be easily modified to achieve this, if desired. The number of buses, lines and units within each zone is also described in Table~\ref{table:size}.
\begin{table}
    \caption{Size of selected instances. \label{table:size}}
    \centering \small
    \begin{tabular}{lrrrrrrrrr}
    \toprule
    & \multicolumn{3}{c}{Total}
    & \multicolumn{3}{c}{Zone 1}
    & \multicolumn{3}{c}{Zone 2} \\
    \cmidrule(l){2-4}
    \cmidrule(l){5-7}
    \cmidrule(l){8-10}
    Instance
    & \hspace{1em} Buses
    & \hspace{1em} Units
    & \hspace{1em} Lines
    & \hspace{1em} Buses
    & \hspace{1em} Units
    & \hspace{1em} Lines
    & \hspace{1em} Buses
    & \hspace{1em} Units
    & \hspace{1em} Lines
    \\
    \midrule
	\texttt{case1888rte} & 1,888 &   297 & 2,531 & 1,113 & 211 & 1,498 & 784 & 86 & 1,033 \\
	\texttt{case1951rte} & 1,951 &   391 & 2,596 & 1,037 & 119 & 1,415 & 923 & 272 & 1,181 \\
	\texttt{case2848rte} & 2,848 &   547 & 3,776 & 1,481 & 226 & 1,957 & 1375 & 321 & 1,819 \\
	\texttt{case3012wp}  & 3,012 &   502 & 3,572 & 1,637 & 322 & 1,938 & 1,388 & 180 & 1,634 \\
	\texttt{case3375wp}  & 3,374 &   596 & 4,161 & 1,649 & 334 & 2,007 & 1,696 & 262 & 2,154 \\
	\texttt{case6468rte} & 6,468 & 1,295 & 9,000 & 2,896 & 544 & 4,049 & 3,588 & 751 & 4,951 \\
	\texttt{case6515rte} & 6,515 & 1,388 & 9,037 & 3,536 & 800 & 4,831 & 2,994 & 588 & 4,206 \\
    \bottomrule
    \end{tabular}
\end{table}

\subsection{Revised release-and-fix heuristic}
\label{subsec:fixrelease}

We recall that the ADMM procedure described in Subsection~\ref{subsec:admm} is not guaranteed to converge in the presence of binary decision variables, since the resulting optimization problem is no longer convex. When solved through ADMM, decentralized MILPs often present an oscillating behavior, where the values of the binary or integral variables flip back and forth between different discrete values, with little change to objective value, and no progress towards global feasibility. Even if such oscillating behavior is not present in particular instances, it may still not be desirable, for performance reasons alone, to repeatedly solve MIQP subproblems at each ADMM iteration.

Based on these observations, Feizollahi et al \cite{feizollahi2015large} proposed a \emph{release-and-fix procedure}, where the continuous relaxation of the problem is solved first, to obtain an initial lower bound, then the procedure alternates between solving the original MIQP (the \emph{release} cycle), and a restricted MIQP which has some (or all) binary variables fixed to particular values (the \emph{fix} cycle), until either the solution becomes globally feasible, or a time limit is reached.

\begin{figure}[t]
    \caption{Revised fix-and-release procedure for Mixed-Integer ADMM. \label{fig:fix-release}}
    \resizebox{\textwidth}{!}{
    \begin{tikzpicture}[auto]
        \tikzstyle{decision} = [diamond, draw, fill=black!7.5, text width=6em, text centered, node distance=3cm, inner sep=0pt]
        \tikzstyle{block} = [rectangle, draw, fill=black!0, text width=7.5em, text centered, rounded corners, minimum height=4em]
        \tikzstyle{line} = [draw, -latex']
        \tikzstyle{cloud} = [draw, ellipse,fill=black!15, node distance=3cm, minimum height=2em]
        \node [cloud] (start) {Start};
        \node [block, right of=start, node distance=7.5em] (miqp) {Solve MIQP subproblems};
        \node [decision, right of=miqp, node distance=10em] (miqpfeas) {Globally feasible?};
        \node [cloud, above of=miqpfeas, node distance=10em] (stop) {Stop};
        \node [block, below of=miqpfeas, node distance=10em] (updatemiqp) {Update multipliers};
        \node [decision, right of=miqpfeas, node distance=10em] (objchange) {Obj. value stagnated?};
        \node [block, right of=objchange, node distance=12.5em] (qp) {Solve QP subproblems};
        \node [decision, right of=qp, node distance=10em] (feasible) {Globally feasible?};
        \node [block, below of=feasible, node distance=10em] (updateqp) {Update multipliers};
        \node [decision, right of=feasible, node distance=10em] (feaschange) {Infeasibility stagnated?};
        \path [line] (start) -- (miqp);
        \path [line] (miqp) -- (miqpfeas);
        \path [line] (miqpfeas) -- node{Yes} (stop);
        \path [line] (miqpfeas) -- node{No} (objchange);
        \path [line] (objchange) -- node{Yes} (qp);
        \path [line] (objchange) |- node [near start] {No} (updatemiqp);
        \path [line] (updatemiqp) -| (miqp);
        \path [line] (qp) -- (feasible);
        \path [line] (feasible) -- node{No} (feaschange);
        \path [line] (feasible) |- node [near start] {Yes} (stop);
        \path [line] (feaschange) |- node [near start] {No} (updateqp);
        \path [line] (feaschange) -- node [near start] {Yes} ++ (0,13em) -| (miqp);
        \path [line] (updateqp) -| (qp);
    \end{tikzpicture}}
\end{figure}
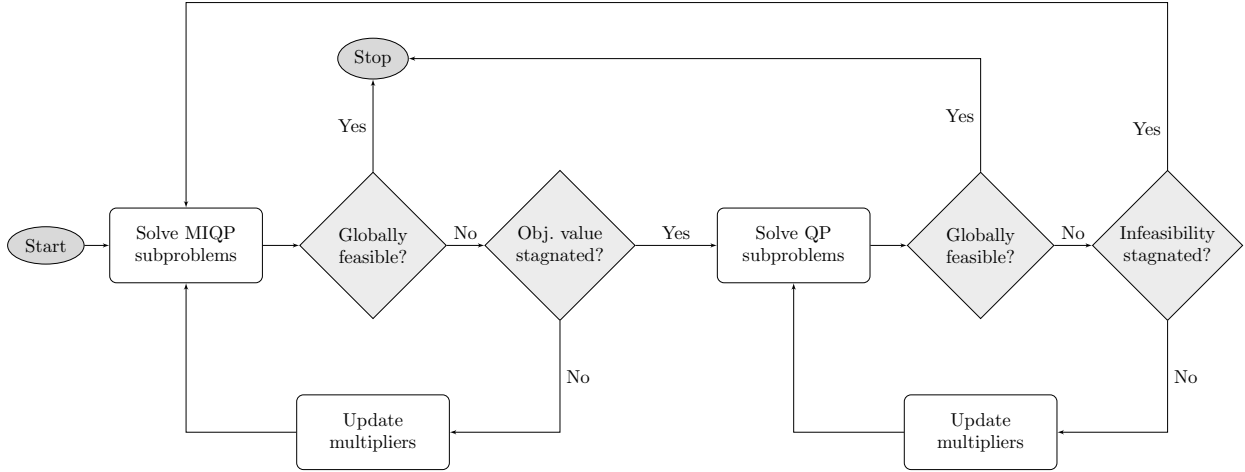

In this subsection, we present a revised version of this \emph{release-and-fix} procedure. The entire process is described in Figure~\ref{fig:fix-release}. At the beginning, we repeatedly solve the original MIQP subproblems and update the ADMM multipliers, according to Algorithm~\ref{alg:admm}. During this \emph{release} cycle, we monitor the changes to the objective value. If, at the end of any ADMM iteration, the global objective value has not changed significantly when compared to the previous iteration, we fix all binary variables to their current optimal values. Then, we switch to the \emph{fix} cycle, where we repeatedly solve restricted these QP subproblems and update the ADMM multipliers. In the \emph{fix} cycle, ADMM is guaranteed to converge to a globally feasible solution, since the problem is convex, unless such a solution does not exist. Assuming that the original instance is feasible, this would indicate that variable fixing performed earlier was not adequate. In this case, the procedure releases all binary variables and returns to the \emph{release} cycle. The procedure repeats until either a globally feasible solution is found, or a time limit is reached. Compared to the method described in \cite{feizollahi2015large}, the revised procedure presented in this subsection mainly differs in the rules used to switch between \emph{release} to \emph{fix} cycles. Here, we propose switches based on changes to objective value and infeasibility, while, in the original method, the switch is made after observing no changes in binary values over the course of 15 iterations, or simply after a fixed number of iterations. The proposed method also differs in its stop criterion. Here, we stop when a globally feasible solution is obtained, whereas in the original method, the search continues even after such a solution is found.

\subsection{Computational Results}

First, we present the computational results for the transmission-constrained version of the problem, which enforces only pre-contingency DC power flow constraints, and no N-1 contingency constraints. Table~\ref{table:tcuc} shows the average running times (in seconds) required to solve the seven instances considered, using either the proposed formulation, or the decentralized phase-angle presented in Subsection~\ref{subsec:dectheta}. For each formulation, the table also shows the primal residual, the number of iterations, and the gap when compared to the optimal solution obtained centrally.

\begin{table}[t]
    \centering \small
    \caption{Benchmark on the Transmission-Constrained Unit Commitment Problem. \label{table:tcuc}}
    \setlength{\tabcolsep}{0.5em}
    \begin{tabular}{lrrrrrrrr}
        \toprule
        & \multicolumn{4}{c}{Proposed}
        & \multicolumn{4}{c}{Phase-Angle}
        \\
        \cmidrule(r){2-5}
        \cmidrule(r){6-9}
        Instance
        &\footnotesize Time (s)
        &\footnotesize Infeas.
        &\footnotesize Iter.
        &\footnotesize Gap (\%)
        &\footnotesize Time (s)
        &\footnotesize Infeas.
        &\footnotesize Iter.
        &\footnotesize Gap (\%)
        \\
        \midrule
        \begin{filecontents*}{tcuc-output.csv}
instance,timeIsf,timeTheta,infeasIsf,infeasTheta,itIsf,itTheta,gapIsf,gapTheta
case1888rte-2z,303.4,3602.7,7.568e-03,5.332e+01,62,300,0.013,nan
case1951rte-2z,204.2,3600.1,9.821e-03,1.603e+01,36,154,0.154,nan
case2848rte-2z,492.4,3600.5,9.369e-03,1.147e+03,47,8,0.299,nan
case3012wp-2z,538.4,3600.4,7.663e-03,5.568e+03,77,7,0.170,nan
case3375wp-2z,474.4,3601.1,8.229e-03,9.527e+03,67,3,0.130,nan
case6468rte-2z,1180.3,3602.3,7.315e-03,2.327e+04,42,3,0.112,nan
case6515rte-2z,1009.7,3605.7,7.462e-03,1.526e+04,45,3,0.008,nan
Average,600.4,3601.8,8.204e-03,7.835e+03,53,68,0.127,nan
\end{filecontents*}
\csvreader[
                head to column names,
                late after line={\\},
                late after last line={\\},
                filter expr={
                    test{\ifnumless{\thecsvinputline}{9}}
                }
            ]{tcuc-output.csv}{}{%
                \texttt{\instance} &
                \timeIsf&
                \infeasIsf&
                \itIsf &
                \gapIsf &
                \timeTheta &
                \infeasTheta &
                \itTheta &
                ---
            }%
          \midrule
          \csvreader[
                  head to column names,
                  late after line={\\},
                  late after last line={\\},
                  filter expr={
                      test{\ifnumgreater{\thecsvinputline}{8}}
                  }
              ]{tcuc-output.csv}{}{%
                \textbf{\instance} &
                \bf \timeIsf&
                \bf \infeasIsf&
                \bf \itIsf &
                \bf \gapIsf &
                \bf \timeTheta &
                \bf \infeasTheta &
                \bf \itTheta &
                \bf ---
              }%
        \bottomrule
    \end{tabular}
\end{table}

Using the proposed formulation, the revised release-and-fix procedure described in Subsection~\ref{subsec:fixrelease} was able to solve all instances well within the 1-hour time limit. In all cases, the procedure ended with a globally feasible solution, with primal residuals within the tolerance. The optimality gaps, when compared to a centralized method, were also relatively small.
On average, the method required 600 seconds, 53 iterations, and produced decentralized solutions which were 0.13\% worse than the optimal central solution. 
Obtaining globally feasible solutions using the decentralized phase-angle formulation, on the other hand, proved very challenging. For all instances, the method exceeded the 1-hour time limit, and terminated without producing any globally feasible solutions. Even for the smallest instances, the primal residuals were still significantly high, indicating that the partial solutions were not compatible enough for the production schedule to be implementable.
With the phase-angle formulation, CPLEX also faced several numerical issues when solving the associated QPs. To mitigate these problems, we used relaxed convergence tolerances for the barrier method (\texttt{CPX\_PARAM\_BAREPCOMP} was set to $10^{-3}$). No such issues were present with the proposed formulation.

\begin{figure}[t]
    \begin{center}
        \caption{Infeasibility over time (instance \texttt{case1888rte-2z}). \label{fig:convergence-a}}
        \includegraphics[width=0.8\textwidth]{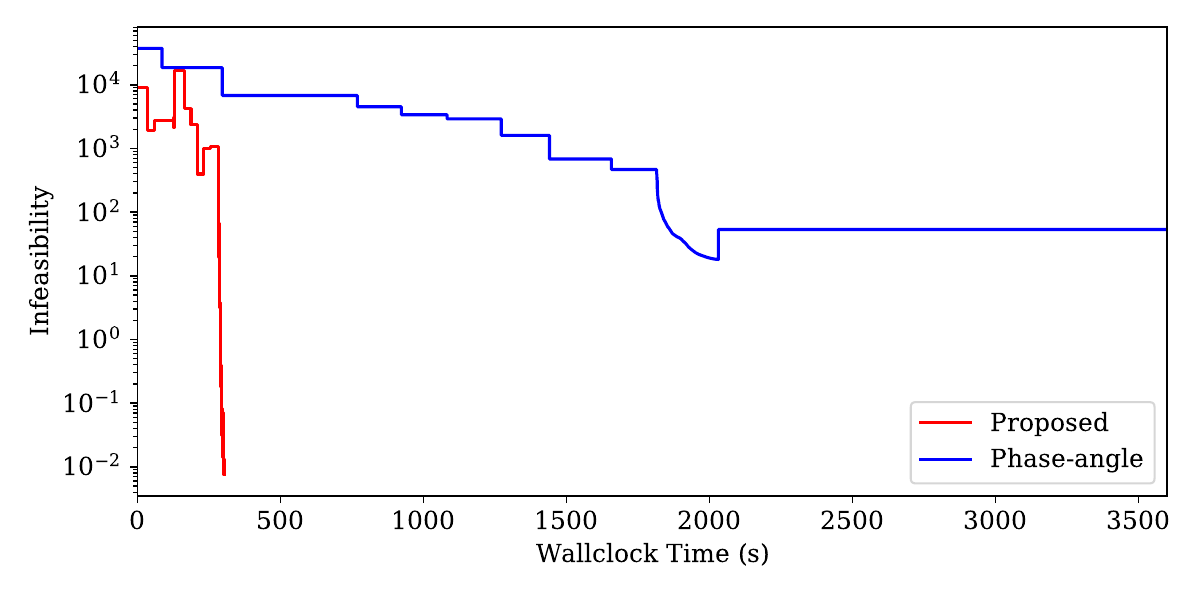}
    \end{center}
\end{figure}

\begin{figure}[t]
    \begin{center}
        \caption{Infeasibility over time (instance \texttt{case6515rte-2z}). \label{fig:convergence-b}}
        \includegraphics[width=0.8\textwidth]{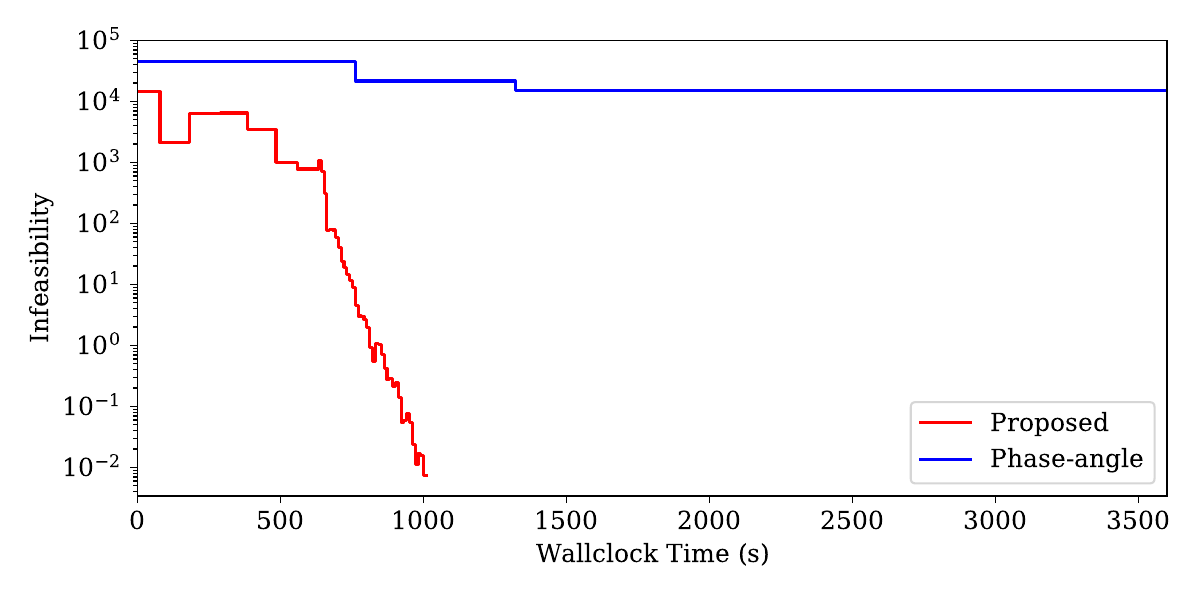}
    \end{center}
\end{figure}

Figures~\ref{fig:convergence-a} and \ref{fig:convergence-b} show the progress towards primal feasibility over time, for instances \texttt{case1888rte-2z} and \texttt{case6515rte-2z}, for both formulations. The longer steps in the chart correspond to the MIQP iterations, whereas the smoother areas correspond to the QP iterations. 
From these figures, it is clear that each MIQP iteration takes significantly longer when using the phase-angle formulation. For both instances, solving the first few phase-angle MIQP iterations took longer than performing the entire ADMM procedure using the proposed formulation.

Now we focus our attention to the security-constrained version of the problem, where N-1 security constraints are enforced using the strategy presented in Section~\ref{sec:security}.
Similarly to the previous table, Table~\ref{table:scuc} shows, for each instance, the average running time (in seconds), the primal infeasibility, the number of iterations and the optimality gap when compared to the central optimal solution. Here, we do not show the results for the phase-angle formulation, since this formulation, in its original form, would require one copy of the $\theta$ variables for each N-1 contingency scenario, resulting in an intractable optimization problem.


Enforcing N-1 security constraints, using the proposed solution method, made some instances harder to solve, but the method was still able to find globally feasible solutions for all instances within the 1-hour time limit. On average, the method required 1170 seconds and 52 iterations to produce solutions that were 0.08\% worse than the centrally-obtained ones. There was no indication that enforcing N-1 security constraints caused any increase in the number of ADMM iterations.

\begin{table}[t]
    \centering \small
    \caption{Benchmark on the Security-Constrained Unit Commitment Problem. \label{table:scuc}}
    \setlength{\tabcolsep}{0.5em}
    \begin{tabular}{lrrrrrrrr}
        \toprule
        & \multicolumn{4}{c}{Proposed}
        \\
        \cmidrule{2-5}
        Instance
        &\footnotesize Time (s)
        &\footnotesize Infeas.
        &\footnotesize Iter.
        &\footnotesize Gap (\%)
        \\
        \midrule
        \begin{filecontents*}{scuc-output.csv}
instance,timeIsf,infeasIsf,itIsf,gapIsf
case1888rte-2z,218.9,4.288e-03,55,0.025
case1951rte-2z,282.5,7.374e-03,56,0.064
case2848rte-2z,1216.4,5.571e-03,46,0.268
case3012wp-2z,814.9,9.137e-03,50,0.077
case3375wp-2z,1173.8,6.132e-03,43,0.024
case6468rte-2z,1541.0,7.735e-03,45,0.038
case6515rte-2z,2954.8,9.132e-03,69,0.086
Average,1171.8,7.053e-03,52,0.083
\end{filecontents*}
\csvreader[
                head to column names,
                late after line={\\},
                late after last line={\\},
                filter expr={
                    test{\ifnumless{\thecsvinputline}{9}}
                }
            ]{scuc-output.csv}{}{%
                \texttt{\instance} &
                \timeIsf&
                \infeasIsf&
                \itIsf &
                \gapIsf &
            }%
          \midrule
          \csvreader[
                  head to column names,
                  late after line={\\},
                  late after last line={\\},
                  filter expr={
                      test{\ifnumgreater{\thecsvinputline}{8}}
                  }
              ]{scuc-output.csv}{}{%
                \textbf{\instance} &
                \bf \timeIsf&
                \bf \infeasIsf&
                \bf \itIsf &
                \bf \gapIsf &
              }%
        \bottomrule
    \end{tabular}
\end{table}

\section{Conclusion}

In this paper, we presented a novel formulation of DC power flows that is specially well-suited for decentralized power systems optimization. Assuming that the transmission network can be separated into zones sharing a small number of tie lines, we proved that it is possible to sparsify the traditional \emph{injection shift factor formulation} by adding a small number of auxiliary constraints. The obtained formulation presents a block-diagonal structure, which lends itself naturally to decomposition methods. We also described how to enforce N-1 security constraints without requiring multiple copies of the decision variables. Computational experiments on a large set of realistic Security-Constrained Unit Commitment instances demonstrated that the proposed formulation performs significantly better than the \emph{decentralized phase angle formulation} used in previous studies. Although we only presented computational results for SCUC, we stress that this formulation is applicable to any power systems optimization problems that employ DC power flows. We are currently investigating its application to other families of problems. Another open question is how to optimally subdivide a large transmission network into multiple smaller zones, for improved computational performance.

\ACKNOWLEDGMENT{We gratefully acknowledge use of the Bebop cluster in the Laboratory Computing Resource Center at Argonne National Laboratory. This work was partially supported by \textbf{Laboratory Directed Research and Development} (LDRD) funding from Argonne National Laboratory, provided by the Director, Office of Science, of the U.S. Department of Energy under Contract No. DE-AC02-06CH11357. This work was also partially supported by the U.S. Department of Energy \textbf{Advanced Grid Modeling Program} under Grant DE-OE0000875.}

\begin{APPENDIX}{}
The original benchmark instances considered in our benchmarks had no zonal information. To separate the buses into two zones, an auxiliary Mixed-Integer Linear Optimization (MILP) problem was solved. The problem has one binary variable $x_b$ for each bus $b \in B$, indicating whether $b$ is a boundary bus, and one binary variable $y_l$ for each $l \in L$, indicating whether line $l$ belongs to zone 1.
The MILP tries to minimize the number of boundary buses, while keeping the number of lines in each zone roughly balanced:
\begin{subequations}
    \label{eq:split}
    \begin{align}
        \text{minimize} \;\;
            & \sum_{b \in B} x_b \label{eq:split-a} \\
        \text{subject to} \;\;
            & y_q + (1 - y_m) \leq 1 + x_b
                & \forall b \in B, q \in L(b),l \in L(b), l \neq q, \label{eq:split-a} \\
            & (1 - y_q) + y_m \leq 1 + x_b
                & \forall b \in B, q \in L(b),l \in L(b), l \neq q, \label{eq:split-b} \\
            & \sum_{l \in L} y_l \leq |L| \left(\frac{1}{2} + \eta\right), \label{eq:split-c} \\
            & \sum_{l \in L} (1 - y_l) \leq |L| \left(\frac{1}{2} + \eta\right), \label{eq:split-d} \\
            & x_b \in \{0,1\} & \forall b \in B \\
            & y_l \in \{0,1\} & \forall l \in L
    \end{align}
\end{subequations}
where $L(b)$ indicates the set of lines incident to $b$ and $\eta$ is a constant which controls the unbalance tolerance.
Equations \eqref{eq:split-a} and \eqref{eq:split-b} guarantee that, if $b$ is not a boundary bus, then all transmission lines directly connected to it must belong to the same zone.
Equations \eqref{eq:split-c} and \eqref{eq:split-d} enfore that the zones have roughly the same number of transmission lines.
Since an optimal solution was not strictly required, a large relative MILP gap tolerance was used to solve this problem.
Note that no attempt was made to keep the number of generators in different zones balanced. The auxiliary Problem~\eqref{eq:split} can be easily modified to achieve this, if desired.
This problem can also be used recursively to subdivide each zone into smaller subzones.
\end{APPENDIX}

\bibliographystyle{informs2014} 
\bibliography{papers.bib} 

\begin{thebibliography}{29}
\providecommand{\natexlab}[1]{#1}
\providecommand{\url}[1]{\texttt{#1}}
\providecommand{\urlprefix}{URL }

\bibitem[{Amaldi et~al.(2014)Amaldi, Coniglio, \protect\BIBand{}
  Gualandi}]{amaldi2014coordinated}
Amaldi E, Coniglio S, Gualandi S (2014) Coordinated cutting plane generation
  via multi-objective separation. \emph{Mathematical Programming}
  143(1-2):87--110.

\bibitem[{Baldick \protect\BIBand{} Chatterjee(2014)}]{BALDICK2014319}
Baldick R, Chatterjee D (2014) Coordinated dispatch of regional transmission
  organizations: Theory and example. \emph{Computers \& Operations Research}
  41:319--332, ISSN 0305-0548,
  \urlprefix\url{http://dx.doi.org/https://doi.org/10.1016/j.cor.2012.12.016}.

\bibitem[{Biskas \protect\BIBand{} Bakirtzis(2004)}]{BiskasBakirtzis2004}
Biskas P, Bakirtzis A (2004) Decentralised security constrained dc-opf of
  interconnected power systems. \emph{IEE Proceedings-Generation, Transmission
  and Distribution} 151(6):747--754.

\bibitem[{Bixby(2002)}]{bixby2002solving}
Bixby RE (2002) Solving real-world linear programs: A decade and more of
  progress. \emph{Operations research} 50(1):3--15.

\bibitem[{Boyd et~al.(2011)Boyd, Parikh, Chu, Peleato, Eckstein
  et~al.}]{boyd2011distributed}
Boyd S, Parikh N, Chu E, Peleato B, Eckstein J, et~al. (2011) Distributed
  optimization and statistical learning via the alternating direction method of
  multipliers. \emph{Foundations and Trends{\textregistered} in Machine
  learning} 3(1):1--122.

\bibitem[{Chen et~al.(2016)Chen, Casto, Wang, Wang, Wang, \protect\BIBand{}
  Wan}]{ChenCastoWangWangWangWan2016}
Chen Y, Casto A, Wang F, Wang Q, Wang X, Wan J (2016) Improving large scale
  day-ahead security constrained unit commitment performance. \emph{IEEE
  Transactions on Power Systems} 31(6):4732--4743.

\bibitem[{Chung et~al.(2011)Chung, Kim, \protect\BIBand{} Hur}]{chung2011multi}
Chung K, Kim B, Hur D (2011) Multi-area generation scheduling algorithm with
  regionally distributed optimal power flow using alternating direction method.
  \emph{International Journal of Electrical Power \& Energy Systems}
  33(9):1527--1535.

\bibitem[{Dall'Anese et~al.(2013)Dall'Anese, Zhu, \protect\BIBand{}
  Giannakis}]{dall2013distributed}
Dall'Anese E, Zhu H, Giannakis GB (2013) Distributed optimal power flow for
  smart microgrids. \emph{IEEE Transactions on Smart Grid} 4(3):1464--1475.

\bibitem[{Dey \protect\BIBand{} Molinaro(2018)}]{dey2018theoretical}
Dey SS, Molinaro M (2018) Theoretical challenges towards cutting-plane
  selection. \emph{Mathematical Programming} 170(1):237--266.

\bibitem[{Dey et~al.(2015)Dey, Molinaro, \protect\BIBand{}
  Wang}]{dey2015approximating}
Dey SS, Molinaro M, Wang Q (2015) Approximating polyhedra with sparse
  inequalities. \emph{Mathematical Programming} 154(1-2):329--352.

\bibitem[{Dey et~al.(2018)Dey, Molinaro, \protect\BIBand{}
  Wang}]{dey2018analysis}
Dey SS, Molinaro M, Wang Q (2018) Analysis of sparse cutting planes for sparse
  milps with applications to stochastic milps. \emph{Mathematics of Operations
  Research} 43(1):304--332.

\bibitem[{Doostizadeh et~al.(2016)Doostizadeh, Aminifar, Lesani,
  \protect\BIBand{} Ghasemi}]{doostizadeh2016multi}
Doostizadeh M, Aminifar F, Lesani H, Ghasemi H (2016) Multi-area market
  clearing in wind-integrated interconnected power systems: A fast parallel
  decentralized method. \emph{Energy Conversion and Management} 113:131--142.

\bibitem[{Erseghe(2014)}]{erseghe2014distributed}
Erseghe T (2014) Distributed optimal power flow using admm. \emph{IEEE
  transactions on power systems} 29(5):2370--2380.

\bibitem[{Feizollahi et~al.(2015)Feizollahi, Costley, Ahmed, \protect\BIBand{}
  Grijalva}]{feizollahi2015large}
Feizollahi MJ, Costley M, Ahmed S, Grijalva S (2015) Large-scale decentralized
  unit commitment. \emph{International Journal of Electrical Power \& Energy
  Systems} 73:97--106.

\bibitem[{Guler et~al.(2007)Guler, Gross, \protect\BIBand{}
  Liu}]{guler2007generalized}
Guler T, Gross G, Liu M (2007) Generalized line outage distribution factors.
  \emph{IEEE Transactions on Power Systems} 22(2):879--881.

\bibitem[{{Ji} \protect\BIBand{} {Tong}(2018)}]{multiareainterchange}
{Ji} Y, {Tong} L (2018) Multi-area interchange scheduling under uncertainty.
  \emph{IEEE Transactions on Power Systems} 33(2):1659--1669,
  \urlprefix\url{http://dx.doi.org/10.1109/TPWRS.2017.2727326}.

\bibitem[{Li et~al.(2015)Li, Shahidehpour, Wu, Zeng, Zhang, \protect\BIBand{}
  Zheng}]{li2015decentralized}
Li Z, Shahidehpour M, Wu W, Zeng B, Zhang B, Zheng W (2015) Decentralized
  multiarea robust generation unit and tie-line scheduling under wind power
  uncertainty. \emph{IEEE Transactions on Sustainable Energy} 6(4):1377--1388.

\bibitem[{Loukarakis et~al.(2014)Loukarakis, Bialek, \protect\BIBand{}
  Dent}]{loukarakis2014investigation}
Loukarakis E, Bialek JW, Dent CJ (2014) Investigation of maximum possible opf
  problem decomposition degree for decentralized energy markets. \emph{IEEE
  Transactions on Power Systems} 30(5):2566--2578.

\bibitem[{Loukarakis et~al.(2015)Loukarakis, Dent, \protect\BIBand{}
  Bialek}]{loukarakis2015decentralized}
Loukarakis E, Dent CJ, Bialek JW (2015) Decentralized multi-period economic
  dispatch for real-time flexible demand management. \emph{IEEE Transactions on
  Power Systems} 31(1):672--684.

\bibitem[{Magn{\'u}sson et~al.(2015)Magn{\'u}sson, Weeraddana,
  \protect\BIBand{} Fischione}]{magnusson2015distributed}
Magn{\'u}sson S, Weeraddana PC, Fischione C (2015) A distributed approach for
  the optimal power-flow problem based on admm and sequential convex
  approximations. \emph{IEEE Transactions on Control of Network Systems}
  2(3):238--253.

\bibitem[{Phan \protect\BIBand{} Sun(2014)}]{phan2014minimal}
Phan DT, Sun XA (2014) Minimal impact corrective actions in
  security-constrained optimal power flow via sparsity regularization.
  \emph{IEEE Transactions on Power Systems} 30(4):1947--1956.

\bibitem[{Reid(1982)}]{reid1982sparsity}
Reid JK (1982) A sparsity-exploiting variant of the bartels—golub
  decomposition for linear programming bases. \emph{Mathematical Programming}
  24(1):55--69.

\bibitem[{Sun et~al.(2013)Sun, Phan, \protect\BIBand{} Ghosh}]{sun2013fully}
Sun AX, Phan DT, Ghosh S (2013) Fully decentralized ac optimal power flow
  algorithms. \emph{2013 IEEE Power \& Energy Society General Meeting}, 1--5
  (IEEE).

\bibitem[{Walter(2014)}]{walter2014sparsity}
Walter M (2014) Sparsity of lift-and-project cutting planes. \emph{Operations
  Research Proceedings 2012}, 9--14 (Springer).

\bibitem[{Wang et~al.(2017)Wang, Wang, \protect\BIBand{}
  Wu}]{wang2017distributed}
Wang Y, Wang S, Wu L (2017) Distributed optimization approaches for emerging
  power systems operation: A review. \emph{Electric Power Systems Research}
  144:127--135.

\bibitem[{Wang et~al.(2016)Wang, Wu, \protect\BIBand{} Wang}]{wang2016fully}
Wang Y, Wu L, Wang S (2016) A fully-decentralized consensus-based admm approach
  for dc-opf with demand response. \emph{IEEE Transactions on Smart Grid}
  8(6):2637--2647.

\bibitem[{White \protect\BIBand{} Pike(2011)}]{ISONE2011}
White M, Pike R (2011) Inter-regional interchange scheduling: analysis and
  options. iso new england.

\bibitem[{Xavier \protect\BIBand{} Qiu(2020)}]{UCJL}
Xavier AS, Qiu F (2020) {UnitCommitment.jl: A Julia/JuMP Optimization Package
  for Security-Constrained Unit Commitment (Version 0.1)}.
  \urlprefix\url{http://dx.doi.org/10.5281/zenodo.4275729}.

\bibitem[{Zimmerman et~al.(1997)Zimmerman, Murillo-S{\'a}nchez,
  \protect\BIBand{} Gan}]{matpower}
Zimmerman RD, Murillo-S{\'a}nchez CE, Gan D (1997) Matpower: A matlab power
  system simulation package. \emph{Manual, Power Systems Engineering Research
  Center, Ithaca NY} 1.

\end{thebibliography}
\end{document}